\documentclass[11pt,fleqn]{amsart}

\usepackage{xcolor}
\usepackage[colorlinks=true,linkcolor=black,urlcolor=black,citecolor=black]{hyperref}
\usepackage{a4wide}
\usepackage{euscript,amsmath,amssymb,amsbsy,mathabx,esint}

\usepackage[arrow,matrix,curve]{xy}\usepackage{tikz}\usetikzlibrary{patterns}
\usepackage{array,longtable,multirow,diagbox,graphicx,enumitem,url,hhline,relsize,scalefnt}
\usepackage{float,caption,capt-of}
\newcounter{msct}[section]\renewcommand{\themsct}{\thesection.\arabic{msct}}

\newenvironment{m-proposition}{\vskip3pt\refstepcounter{msct}\trivlist \itemindent0pt\item[\hskip\labelsep\bf Proposition \themsct]\it\ignorespaces}{\endtrivlist\vskip3pt}

\newenvironment{m-corollary}{\vskip3pt\refstepcounter{msct}\trivlist \itemindent 0pt\item[\hskip\labelsep\bf Corollary \themsct]\it\ignorespaces}{\endtrivlist\vskip3pt}

\newenvironment{m-lemma}{\vskip3pt\refstepcounter{msct}\trivlist \itemindent 0pt\item[\hskip\labelsep\bf Lemma \themsct]\it\ignorespaces}{\endtrivlist\vskip3pt}

\newenvironment{m-definition}{\vskip3pt\refstepcounter{msct}\trivlist \itemindent0pt\item[\hskip\labelsep\bf Definition \themsct]\ignorespaces}{\endtrivlist\vskip3pt}

\newenvironment{m-remark}{\vskip3pt\refstepcounter{msct}\trivlist \itemindent0pt\item[\hskip\labelsep\bf Remark \themsct]\ignorespaces}{\endtrivlist\vskip3pt}

\newenvironment{m-thank}{\vskip3pt\trivlist \itemindent0pt\item[\hskip\labelsep\it Acknowledgments]\ignorespaces}{\endtrivlist\vskip3pt}

\newenvironment{m-proof}{\vskip1pt\trivlist \itemindent0pt\item[\hskip\labelsep\it Proof.]\ignorespaces}{\hfill$\Box$\endtrivlist\vskip3pt}

\newcounter{meqn}[section]\renewcommand{\themeqn}{\thesection.\arabic{meqn}}
\newenvironment{m-eqn}[1]{\vskip3pt\refstepcounter{meqn}\trivlist\itemindent0pt%
\item[]\ignorespaces\hfill$\displaystyle #1$\hfill\hbox{\rm(\themeqn)}}{\endtrivlist\vskip3pt}

\newcommand{\bibauth}[2]{\textrm{{#2}~{#1}},}
\newcommand{\bibtitl}[1]{``\textrm{#1},"}
\newcommand{\bibjnyp}[4]{\textrm{#1} {#2}, {pp. #4}, {#3,}}

\newcommand{\bibbook}[4]{\textit{#1}, {#2}, {#3}, {#4}.}

\numberwithin{equation}{section}

\let\mt\mapsto

\let\blt\bullet

\let\disp\displaystyle
\let\aph\alpha

\let\gma\gamma
\let\ges\geqslant   \let\les\leqslant
\let\kp\kappa
\let\lda\lambda

\let\nit\noindent

\let\veps\varepsilon

\newcommand{\exc}{{\rm ex}}

\newcommand{\CC}{C}
\newcommand{\cp}{c_+}\newcommand{\Cp}{C_+}

\newcommand{\aaa}{a}

\newcommand\bone{{1\kern-0.57ex\rm l}}
\newcommand{\bl}{{\rm bl}}
\newcommand{\bvp}{{\rm BVP}}\newcommand{\ivp}{{\rm IVP}}

\newcommand{\gss}{{\rm est}}

\newcommand{\mxerr}{{\rm maxerr}}
\newcommand{\ort}{\mathrel{{\vrule width5pt height0.5pt depth0pt\vrule width0.5pt height7pt depth0pt\,}}}

\newcommand{\ouset}[3]{\underset{#1}{\overset{#2}{#3}}}
\newcommand{\rd}{{\rm d}}

\newcommand{\nn}{n}
\newcommand{\pp}{p}\newcommand{\PP}{B}
\newcommand{\rr}{r}

\newcommand{\ww}{w}
\newcommand{\xx}{x}
\newcommand{\yy}{y}\newcommand{\YY}{Y}\newcommand{\yb}{y_\bullet}
\newcommand{\yp}{y_{+}}\newcommand{\ym}{y_{-}}
\newcommand{\yip}{y_{+}}\newcommand{\yim}{y_{-}}
\newcommand{\zz}{z}\newcommand{\ZZ}{Z}

\begin{document}

\title[On a {\bvp} related to heat-radiation]{On a boundary value problem \\ related to radiative heat transfer} 
\author{Mihai Halic}
%\email{mihai.halic@gmail.com}\address{EMR-matrix, 310029 Arad, Romania}

\begin{abstract} 
We consider the boundary value problem $\,y''=\kp^2y^n, y(0)=1,y(1)=0,\,$ relevant to radiative heat transfer, chemical reactions, etc. We use upper/lower envelopes to study the qualitative and quantitative dependence of its solution on parameters. The envelopes allow estimating the initial derivative, so that the problem can be transformed into an easily solvable initial value problem. Moreover, the knowledge of the approximate value of $\yy'(0)$ allows constructing an explicit approximating function for $\yy(x)$.
\end{abstract}

\thispagestyle{empty}\setcounter{page}{0}
\newpage 
\maketitle
\pagestyle{myheadings}
\markboth{\sc On a {\bvp} related to heat-radiation}{\sc On a {\bvp} related to heat-radiation}

\section{Introduction}

We consider the boundary value problem (\bvp), 
\begin{m-eqn}{
\yy''=\kp^2\cdot\yy^\nn,\,\yy(0)=1,\;\yy(1)=0,
\;\text{($\nn>1, \kp>0$ are real parameters)}. 
}\label{eq:bvp}
\end{m-eqn} 
It is rooted in Stefan-Boltzmann's law~\cite[\S21.2, (21.5)--(21.8)]{mod}: 
\begin{m-eqn}{
\frac{\rd^2\yy}{\rd\xx^2}=\frac{{1}}{N}(\yy^4-g), \quad 
\begin{array}{lr}
N=\text{conduction-to-radiation parameter},\\ 
g=O((\text{initial temperature})^{-4}).
\end{array}
}\label{eq:sb}
\end{m-eqn}  
Recall~[\textit{id.}, pp.\,1-2]: ``\textrm{radiative heat transfer rates are generally proportional to differences in temperature to the fourth (\emph{or higher}) power. Therefore, radiative heat transfer becomes more important with rising temperature levels and may be totally dominant.}" In this situation, $g, N$ are small, leading to~\eqref{eq:bvp} for $g\to0$. Numerical computations require \textit{guessing first} the temperature profile, due to the non-linearity of the equation [\textit{ib.}, pp.\,777,\,779]. 
The arbitrary exponent $\nn$, not necessarily $4$ or integer, is justified by investigations~\cite{ch+pa,len+men,mppt} of the blackbody radiation problem for  nonextensive systems --obeying Tsallis rather than Boltzmann-Gibbs statistics-- possessing large or infinite numbers of degrees of freedom. Super-Planckian thermal radiation occurs also~\cite{lvc,zzw} when restricted to narrow frequency bands.

The equation appears in other contexts as well: can be viewed as the $1$-dimensional, stationary, generalized Klein-Gordon equation \cite{nrt,keg}; governs the concentration of a substance disappearing in an $n^{\rm th}$-order isothermal reaction~\cite{gav,me-ar}. The use of hypergeometric functions in the latter reference gives limited insight into the behaviour of the solution. It's surprising that, apparently, combined analytical/computational methods are rarely used.

For the author, a further motivation to study~\eqref{eq:bvp} comes from the fact that he analysed~\cite{hal-aip} the equation~\eqref{eq:sb} and derived an explicit approximate analytical solution, which involves no guessing. The precision of the estimate is satisfactory (cf. [\textit{id.}, Table 1]) as long as $g$ is not too small; however, the error increases for $g\to0$. 

Although there are numerous algorithms for solving {\bvp}s in general (see~\cite{amr}), boundary layer problems challenging because a single computation or plot may require (on usual computers) hours or more, especially for non-expert programmers. This inconvenience disappears if one has an explicit expression at hand. 
Direct numerical approaches to~\eqref{eq:bvp} are challenging for large $\kp, \nn$: the solution varies very steeply near $\xx=0$. Here we determine  explicit upper/lower envelopes which squeeze in between the exact solution and clarify the dependence on parameters. At the boundary layer, where the steep variation takes place, they yield accurate values for the exact solution. Also, allow estimating the (unknown) initial derivative $\yy'(0)$, thus transforming the {\bvp} into a routine {\ivp} which can be numerically solved on a usual laptop, with standard software. The knowledge of the approximate value of the initial derivative, in its turn, yields an explicit function which accurately approximates $\yy_\exc$ along the boundary layer interval, removing the necessity of numerical solvers.

The article is structured as follows: in \S\ref{sct:theory} we set up the necessary analytical package;  second, in \S\ref{sct:heat}, we apply them to~\eqref{eq:bvp}, construct envelopes, deduce the approximate value of $\yy'(0)$, and construct the explicit approximating function. The code is shared\footnote{The MAPLE code used by the author for numerical computations is available at: \\ \url{https://drive.google.com/file/d/1Gdd7uoNC9bELjI0uYysYmklvXjinKvQK/view?usp=sharing}} so that one may analyse its own favourite situation.

%%%%%%%%%%
%%%%%%%%%%

\section{Analytical toolbox}\label{sct:theory}

Numerical algorithms are highly effective as long as the  solution doesn't vary too dramatically. In such cases, to avoid errors, tiny meshes require high computer power. A numerical solution is exact if it differs from the analytic one by a small amount. We declare that \emph{any function sufficiently close} to the analytical solution is \emph{numerically exact}.  

\subsection{Criteria for upper/lower envelopes} 

We consider equations of the form: 
\begin{m-eqn}{
\begin{array}{l | l | l}
\underline{\text{BVP-2}}
&\underline{\text{BVP-1}}
&
\\ 
\yy''=f(\xx,\yy)
&\yy'=\sqrt{2}\cdot F(\xx,\yy)
&\text{$f, F$ are differentiable} 
\\ 
\yy(\xi_0), \yy(\xi_1)\;\text{given},
&\yy(\xi_0), \yy(\xi_1)\;\text{given}.
&\text{on}\;[\xi_0,\xi_1]\times[\Xi_0,\Xi_1].
\end{array}
}\label{eq:BVP12}
\end{m-eqn}
Let $\yy_\exc$ be the (analytically-exact) solution. Since BVP-1 is of first-order with two boundary conditions, $F$ typically involves a parameter to ensure the problem is solvable. Our approach to solving them numerically is to replace $f, F$ by functions $g, G$, respectively, such that: 
\begin{itemize}
\item 
$|g-f|,\;|G-F|$ are small. (The solutions of the modified equations approximate $y_\exc$.)
\item 
The equations $\;\yy''=g(\xx,\yy)$ and $\yy'=\sqrt{2}\cdot G(\xx,\yy)\;$ are explicitly integrable.
\end{itemize}
The first requirement is loose, the second is restrictive. In fact, to estimate the accuracy of the approximation, we need \emph{two functions} $g_\pm$ (resp. $G_\pm$), so that the solutions $\yy_\pm$ satisfy:
\begin{m-eqn}{
\ym\les y_\exc \les \yp.
}\label{eq:ypm}
\end{m-eqn}
The error is $\,\mxerr:=\disp\max_{x\in[a,b]}|\yp(x)-\ym(x)|.$ Below is the tool for checking~\eqref{eq:ypm}. 

\begin{m-proposition}\label{prop:+-} 
Let $f,g$ be continuous on $[\xi_0,\xi_1]\times[\Xi_0,\Xi_1]$ and $y, z$ solve respectively:
\vskip.5ex\centerline{
$
(\bvp_{\!f})\, \yy''=f(\xx,\yy),\;\yy(\xi_0), \yy(\xi_1)\,given;
\;
(\bvp_{\!g})\, \zz''=g(\xx,\zz),\;\zz(\xi_0), \zz(\xi_1)\,given.
$
}\vskip.5ex
\nit Then $\yy(x)\les \zz(x), \forall\xx\in\![\xi_0,\xi_1]$, if the following conditions are satisfied: 

\nit{\rm(a)} $f(\xx,\rr)\!\ges\! g(\xx,\rr), \forall\rr\in\![\xi_0,\xi_1],$ and either $f$ or $g$ is strictly increasing in $\rr$;

\nit{\rm(b)} $y(\xi_0)\les z(\xi_0),\;y(\xi_1)\les z(\xi_1).$ 
\end{m-proposition}

\begin{m-proof} If $\{\xx\mid \yy(\xx)>\zz(\xx)\}\neq\emptyset$, contains a maximal interval $(c,d)\subset[\xi_0,\xi_1]$; so we have $\yy(c)=\zz(c)$, $\yy(d)=\zz(d)$. We claim: $\,\yy''(\xx)>\zz''(\xx),\,\forall\,\xx\in(c,d).$ Suppose $g$ is increasing in the $2^{\rm nd}$ variable: $\,\yy''(\xx)=f(\xx,\yy(\xx))\ges g(\xx,\yy(\xx))>g(\xx,\zz(\xx))=\zz''(\xx)$, similarly for $f$. 

We deduce:  
$0=\int_c^d(\yy'(\xx)-\zz'(\xx))\rd\xx=(d-c)\cdot(\yy'(d)-\zz'(d))-\int_c^d(\xx-c)(\yy''(\xx)-\zz''(\xx))\rd\xx,$ 
so $\yy'(d)-\zz'(d)>0$; the same is true on some $(d-\veps,d+\veps)$. But the intermediate value theorem for $\yy-\zz$ implies $\yy'(\tau)-\zz'(\tau)<0$ for some $\tau\in(d-\veps,d)$, a contradiction. \end{m-proof}

%%%%%%%%%%
%%%%%%%%%%

\subsection{Autonomous BVP} 

This means that $f=f(\yy)$, with $f$ assumed to be differentiable. A linear change of variables leads to 
$\;\yy''=f(\yy),\; \yy(0)=0,\;\yy(1)=1.$ 
Let $F$ be the anti-derivative of $f$, with $F(0)=0$, so $\yy_\exc$ satisfies the \ivp:  
\begin{m-eqn}{
\yy'=\sqrt{2}\cdot\sqrt{F(\yy)+\CC},\,\yy(1)=1,
}\label{eq:C}
\end{m-eqn} 
where $\CC=\CC_\exc$ is determined (uniquely, in general) by $\yy(0)=0$. Note that, if $F\ges0$ and $\disp\int_0^1\frac{\rd\yy}{\sqrt{F(\yy)}}\ges\sqrt{2}$, then the equation admits a unique solution: the integral is strictly decreasing with $\CC$,  it's small for $\CC\gg0$, so equality is achieved for only one value $\CC>0$. 

\subsubsection{Boundary layer (bl)}\label{sssct:bd-layer} 

It reflects that $\yy=\yy(\xx)$ varies steeply near a vertical line, like an asymptote; we consider $\xx=1$ and $\yy=\yy_\exc$. 
The condition $\yy'_\exc(1)\gg0$ is not good enough, the graph may `move away' quickly.  Instead, we require: 
\begin{m-eqn}{
\yy'_\exc(\yy^{-1}_\exc(0.9))>50\approx\tan(89^\circ)
\;\Leftrightarrow\;(\yy^{-1}_\exc)'(0.9)<0.02\approx\tan(1^\circ). 
}\label{eq:bd-yexc} 
\medskip\newline By inserting into~\eqref{eq:C}, $\yy_\exc$ has boundary layer if $\,F(0.9)\ges1250.$
\end{m-eqn}

For $\xi=\yy^{-1}_\exc(0.9)$, we obtain: $\disp\frac{1-0.9}{1-\xi}>50\Rightarrow1-\xi<0.002.$ Thus the quantitative meaning is that the function drops $10\%$ of its total variation (from $1$ to $0$) on an interval of length less than $0.2\%$ of the length of its domain. (Note: the choices 10\% and 0.2\% can be adjusted according to the situation at hand.) Such functions have $\ort$-shaped graphs (see Fig.~\ref{fig:ortshape}): stretch near the $X$-axis, then suddenly raise.

%%%%%%%%%%
%%%%%%%%%%

\subsubsection{Upper/lower envelopes}%\label{sssct:y+-} 

We assume that $f$ is positive on $[0,1]$; thus $F\ges0$, too. Note that, in this case, $\yy_\exc$ is convex and increasing. 

\begin{m-corollary}\label{cor:y+-}
\begin{enumerate}
\item[\rm(i)] 
The solution $\yp$ of~\eqref{eq:C} with $\CC=0$ is an upper envelope of $\yy_\exc$. It satisfies $\yp\ges\yy_\exc$ on $[0,1]$. 
\item[\rm(ii)] 
One has $\CC_\exc<\Cp$, where $\Cp$ is defined by~\eqref{eq:h}. The solution $\ym$ of~\eqref{eq:C} with $\CC=\CC_+$ is a (partial) lower envelope of $\yy_\exc$; satisfies $\ym\les\yy_\exc$ on its domain of definition.
\item[\rm(iii)] Suppose $f$ is positive, increasing. Then we have $\CC_-:=F(\yp(0))<\CC_\exc$.
\end{enumerate}
\end{m-corollary}

\begin{m-proof}
(i) This follows from Prop.~\ref{prop:+-}, since $\sqrt{F(\yy)}<\sqrt{F(\yy)+\CC}$. The solution $\yp$ is defined on $[0,1]$, because it's increasing, greater than $\yy_\exc$, so it remains bounded. 

\nit(ii) The upper envelope determines an upper bound $\Cp$ for $\CC$. Prop.~\ref{prop:+-} implies that $\ym\les \yy_\exc$.\\[1ex] 
\begin{tabular}{l}
\begin{minipage}[c]{.57\textwidth}
 The picture shows that $\yy'_\exc(0)$ is smaller than the slope of the secant joining $(0,0)$ and $(\xi,\yp(\xi))$.  
\begin{m-eqn}{\kern-1ex
\sqrt{2\CC_\exc}\les\frac{\yp(\xi)}{\xi}
\,\Rightarrow\,
\CC_\exc\les \Cp\!:=\!\frac{1}{2}\Bigl(\frac{\yp(\xi)}{\xi}\Bigr)^2\!.
}\label{eq:h}
\end{m-eqn}
This holds true for any $\xi\in(0,1)$. One is interested in the `optimal' $\xi$, for which the ratio is the smallest (see Remark~\ref{rk:C+}). 
\end{minipage} 
\begin{minipage}[c]{.4\textwidth}
\begin{m-eqn}{
\scalebox{0.85}{
\begin{tikzpicture}[baseline=(baseline point)]
\coordinate (baseline point) at (0, 1);
\draw[thick,->] (-1,0) -- (3.75,0) node[below] {$x$};
\draw[thick,->] (0,-.5) -- (0,2.5) node[left] {$y$};
\draw[thick,domain=0:3.67] plot (\x,{0.01*\x*\x*\x+0.2*\x});
\node at (3.25,0.5) {$y = \yy_{\exc}$};
\draw[thick,domain=0:3.67] plot (\x,{0.01*\x*\x*\x+0.1*\x+1});
\node at (1,1.35) {$\yy = \yp$};
\node at (-0.5,-0.25) {$(0,0)$}; \node at (0,0) {$\bullet$};
\node at (-0.85,1) {$(0,\yp(0))$};
\node at (2.5,-0.33) {$\xi$};
\draw[very thick,dashed] (2.5,0) -- (2.5,0.01*2.5*2.5*2.5+0.1*2.5+1) node[above] {$\yp(\xi)$};
\draw[very thick,dashed] (0,0) -- (2.5,0.01*2.5*2.5*2.5+0.1*2.5+1);
\node at (2.5,0) {$\bullet$};\node at  (2.5,0.01*2.5*2.5*2.5+0.1*2.5+1)  {$\bullet$};
\end{tikzpicture}
}}\captionof{figure}{$\CC_\exc\les\CC_+$}\label{fig:Cp}
\end{m-eqn}
\end{minipage}
\end{tabular}

\nit(iii) We have $(\yy_\exc'-\yp')'(x)= f(\yy_\exc(\xx))-f(\yp(\xx))\les0$, so $\yy_\exc'-\yp'$ is decreasing. Since $(\yy_\exc'-\yp')(1)>0$, the same holds true for $\xx=0$, hence $\sqrt{2\CC_\exc}>\sqrt{2F(\yp(0))}$. 
\end{m-proof}

\nit\begin{minipage}{.98\textwidth}
\begin{minipage}[c]{.62\textwidth}
\begin{m-remark}\label{rk:C+}
The smallest $\Cp$ in~\eqref{eq:h} corresponds to the `optimal' $\xi$, which satisfies the equation:  
\\[.5ex]\centerline{
$\;2\xi^2\cdot F(\yp(\xi))=\yp(\xi)^2\;\;\Rightarrow\;\;\Cp= F(\yp(\xi)).$
}\\[.5ex]  
In the bl-case, $\xi\approx1$ due to the $\ort$-shape, so $\xi$ (approximately) solves $2F(y)=y^2$. 
\end{m-remark}
\end{minipage}
\begin{minipage}[c]{.33\textwidth}
\centering
\textscale{.73}{ 
\includegraphics[height=0.8\textwidth]{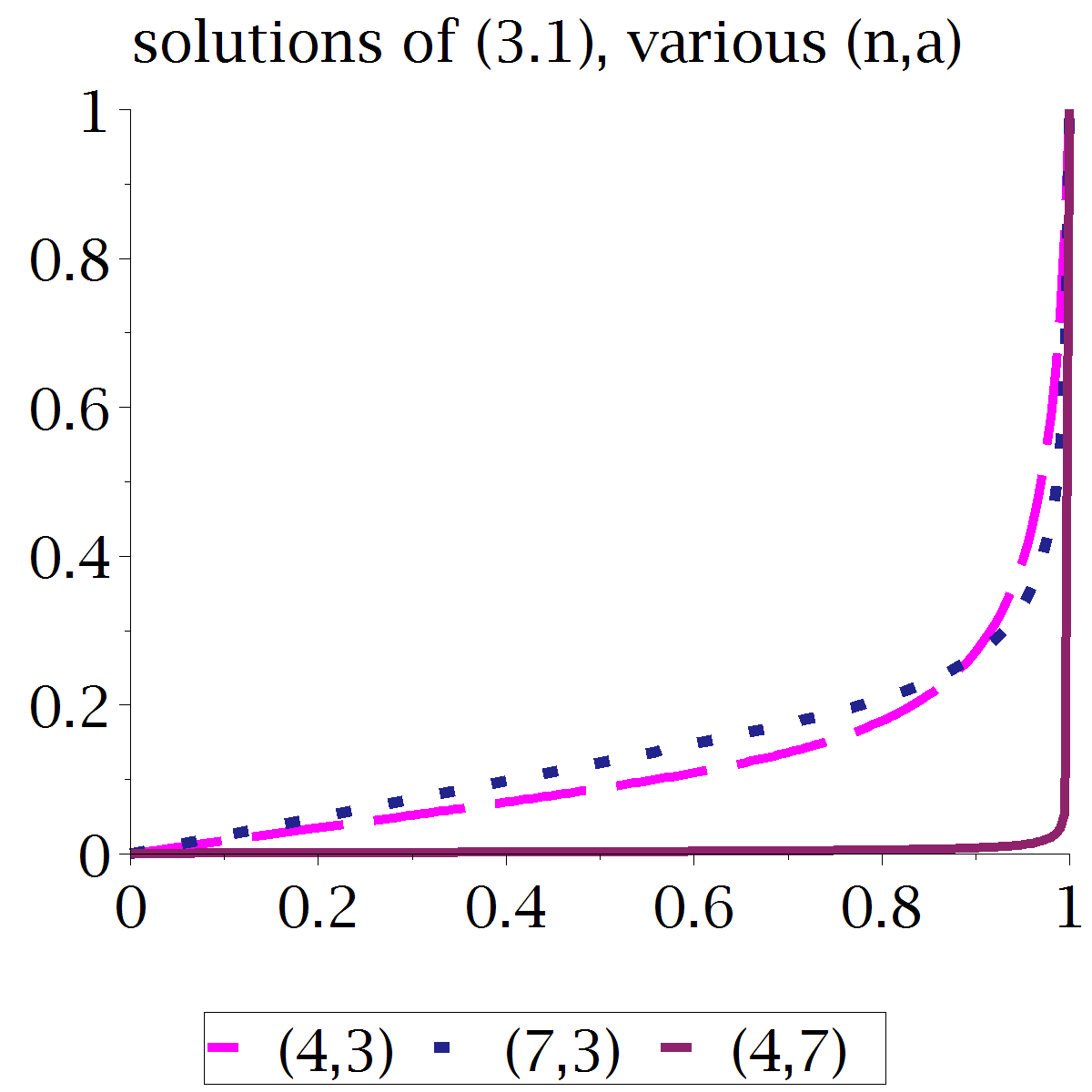}
\captionof{figure}{{sols.\,of\,\eqref{eq:h2} }}\label{fig:ortshape}
}
\end{minipage}
\end{minipage}

%%%%%%%%%%
%%%%%%%%%%

\subsubsection{Numerical point of view}\label{sssct:RK}

The envelope technique is strengthened when enhanced with computer-power. The solution of~\eqref{eq:BVP12} satisfies~\eqref{eq:C} with an unknown $\CC\in[0,\Cp]$. We turn around this standpoint. \emph{For any} such $\CC$, the solution $\yy_\CC$ of~\eqref{eq:C} satisfies $\yy''=f(\yy)$, so we have: 
\begin{m-eqn}{\text{\it 
$\yy_\CC$ is an upper (resp. lower) envelope of $\yy_\exc\,\Leftrightarrow\,\yy_C(0)\ges0$ (resp. $\les0$).
}}\label{eq:yC}
\end{m-eqn}
The benefit of this viewpoint  is that numerical algorithms for {\ivp}s ---we use Runge-Kutta (RK)--- are more robust than for {\bvp}s. It allows verifying~\eqref{eq:yC} \emph{if we pick the correct value of $\CC$!} For this reason, having (tight) a priori bounds $\CC_\pm$ are essential. 
Estimating the difference (error) between $\yy_\CC$ and $\yy_\exc$ is easy: the smaller $|\yy_\CC(0)|$, the more accurate the estimate, uniformly on $[0,1]$: $\disp\;\mxerr=\max_{x\in[0,1]}|\yy_\CC(x)-\yy_\exc(x)|=|\yy_\CC(0)|.$  

\nit Indeed, both $\yy_\CC, \yy_\exc$ satisfy $\yy''=\kp^2 \yy^\nn$. If $\yy_\CC(0)>0$, then $\yy_\CC$ is greater than $\yy_\exc$,  so $\yy_\CC''-\yy_\exc''>0$, hence $\yy_\CC'-\yy_\exc'$ is increasing. But $\yy_\CC\ges \yy_\exc$ implies $\yy_\CC'(1)-\yy_\exc'(1)<0$, so $\yy_\CC'-\yy_\exc'<0$ on $[0,1]$, thus $\yy_\CC-\yy_\exc$ is decreasing.

%%%%%%%%%%
%%%%%%%%%%

\section{Application: Stefan-Boltzmann-type power-law}\label{sct:heat}

We analyse the dependence of the solution of 
\begin{m-eqn}{
\yy''=\kp^2\cdot \yy^\nn,\quad \yy(0)=0,\;\yy(1)=1,\;\;(\nn,\kp>0), 
}\label{eq:h2}
\end{m-eqn} 
on the parameters and compute it accurately, especially at the boundary layer. (We substituted $\xx\mt 1-\xx$ in~\eqref{eq:bvp}.) It's equivalent to the {\ivp} below, where the unknown $\CC=\CC_\exc>0$ is determined by the condition $\yy(0)=0$:
\begin{m-eqn}{
\yy'=\kp\sqrt{2}\cdot \sqrt{\frac{\yy^{\nn+1}}{\nn+1}+\CC},\;\;\yy(1)=1,
}\label{eq:h1}
\end{m-eqn} 
Our approach involves the following steps: 
\begin{enumerate}
\item[\rm(i)] Construct explicit upper and lower envelopes for the exact solution $\yy_\exc$. Using them, we geometrically define the boundary layer interval and determine upper/lower bounds for the initial derivative $\CC_\exc$.
\item[\rm(ii)] Use numerical data to deduce the approximate value of $\CC_\exc$ in terms of $\nn, \kp$. Thus the {\bvp}~\eqref{eq:h2} is transformed into the {\ivp}~\eqref{eq:h1}, which is much easier to solve numerically. 
\item[\rm(iii)] Use the approximate value of $\CC_\exc$ to construct an explicit function which accurately approximates $\yy_\exc$ on the bl-interval, removing the necessity of numerical solvers. 
\end{enumerate}

\subsection{Analytical consideration} 

We are going to construct upper/lower envelopes for $\yy_\exc$.

\begin{m-proposition}\label{prop:estim-C}
\begin{enumerate}
\item[\rm(i)] 
The solution of~\eqref{eq:h2} has boundary layer for 
\vskip.5ex\centerline{$\disp\kp\ges\kp_\bl:=36\sqrt{\frac{\nn+1}{0.9^{\nn+1}}}.$}\vskip.5ex  
We call  $\aph_\bl:=\log_\nn(\kp_\bl)$ \emph{the boundary layer exponent}.
\item[\rm(ii)] 
The function 
$\;\yip(x)\!:=\!{\Bigl[ 1\!+\!\frac{\kp(\nn-1)}{\sqrt{2(\nn+1)}}(1-\xx) \Bigr]}^{-\frac{2}{\nn-1}}\!$ 
is an upper envelope of $y_\exc$ on $[0,1]$.

\item[\rm(iii)] Let $\disp\xi:=\frac{\nn-1}{\nn+1}\les1,$ $\disp\Cp:=\frac{1}{2}\bigl[ \frac{\yip(\xi)}{\kp\xi} \bigr]^2
=\frac{1}{2\kp^2\xi^2}
\Bigl[1+\frac{\sqrt{2}(\nn-1)\kp}{(\nn+1)^{3/2}}\Bigr]^{-\frac{4}{\nn-1}},$ and 
\item[] $\disp\cp:=[(\nn+1)\Cp]^{\frac{1}{\nn+1}}$. We define   
\vskip.5ex\centerline{
$\disp\yb(x):={\biggl[ (1+\cp)^{-\frac{\nn-1}{2}}+\frac{\kp(\nn-1)}{\sqrt{2(\nn+1)}}(1-\xx) \biggr]}^{-\frac{2}{\nn-1}}{\kern-1ex} -\cp.$ 
}\vskip.5ex  
\nit For $\disp\lda:=\frac{\yip(0)}{\yip(0)-\yb(0)}$, the function 
$\yim:=(1-\lda)\yip+\lda \yb$ satisfies $\yim(0)=0$, $\yim(1)=1$. Moreover, $\yim$ has the following properties: 
\begin{itemize}
\item it is increasing, convex, positive, $\yip-\yim$ is decreasing on $[0,1]$;
\item it's a lower envelope of $y_\exc$.
\end{itemize}
\item[\rm(iv)] Let 
$\disp\CC_-:=\frac{\yip(0)^{\nn+1}}{\nn+1}
=\frac{1}{\nn+1}{\Bigl[ 1+\frac{\kp(\nn-1)}{\sqrt{2(\nn+1)}} \Bigr]}^{-\frac{2(\nn+1)}{\nn-1}}.$ 
Then it holds true:\\\centerline{$\CC_-\les\CC_\exc\les\CC_+$.}
\end{enumerate}
\end{m-proposition}

\begin{m-proof}
The claim~(i) follows from~\eqref{eq:bd-yexc}. Cor.~\ref{cor:y+-} yields~(ii). 

\nit(iii) The optimal value of $\xi$ in~\eqref{eq:h} is $\frac{\nn-1}{\nn+1}\!+\!\frac{\sqrt{2}}{\kp\sqrt{\nn+1}}$, but we choose $\frac{\nn-1}{\nn+1}$ for convenience. Since 
\\[.5ex]\centerline{
$\disp\sqrt{2}\kp\sqrt{\frac{y^{\nn+1}}{\nn+1}+\Cp}
\les\frac{\sqrt{2}\kp}{\sqrt{\nn+1}}\cdot\bigl( \yy+\cp \bigr)^{(\nn+1)/2},$
}\\[.5ex]  
Cor.~\ref{cor:y+-} implies that the solution $\yb$ of 
$\,\yy'=\frac{\sqrt{2}\kp}{\sqrt{\nn+1}}\cdot\bigl( \yy+\cp \bigr)^{(\nn+1)/2},\;\yy(1)=1,$ 
is a lower envelope of $\yy_\exc$. Prop.~\ref{prop:+-} implies $\yim\les y_\exc$: indeed, $\yim$ is increasing, thus positive, and 
\\[.25ex]\centerline{$\;\yim''(\xx)=\kp^2\big((1-\lda)\yip(\xx)^\nn+\lda (\yb(\xx)+\cp)^\nn\big)\ouset{ineq.}{Jensen}{>} \yim(\xx)^\nn.$}\\[.25ex] 
The functions $\yip, \yim$ are explicit, one verifies that the difference of their derivatives doesn't vanish; thus it's negative, because so is at $\xx=1$.

\nit(iv) This follows from Cor.~\ref{cor:y+-}(iii), because $\xx\mt\xx^\nn$ is increasing. 
\end{m-proof}

\nit The result is not only theoretical. It yields accurate values of $\yy_\exc$ at the boundary layer ---it's squeezed between $\yim$ and $\yip$---, for parameters which are troublesome for `pure force' numerical approaches (see Table~\ref{tab:bdry}). The a priori estimates for $\CC_\exc$ are the key for successfully running the `shooting method'. Without, numerical tests fail, computer stalls. 

%%%%%%%%%%
%%%%%%%%%%

\subsection{Numerical analysis}\label{sct:num-k}

We confront the theoretical results against numerical data: 
\\\nit -- verify the explicit envelopes constructed above; 
\\\nit-- improve them (see \S\ref{sssct:empir}) using the numerical approach described in~\S\ref{sssct:RK}.  

\subsubsection{Using envelopes} 

They are useful to compute explicit, accurate values at the boundary layer, where $\yy_\exc$ varies very steeply. We start by geometrically defining the region in the XY-plane corresponding to the boundary layer. 
\\[1ex]\begin{tabular}{l}
\begin{minipage}[c]{.5\textwidth} 
Let $\PP$ be the intersection between the X-axis and the tangent line to $\yim$ at $x=1$, so $$1-\PP=(\yim'(1))^{-1}.$$ 
The triangle formed with the line $x=1$ approximates the boundary component of $\yy_\exc$. Let: 
$\;\pp=\log_{10}(\yim'(1)),\;\;\text{so}\;\;\PP=1-10^{-\pp}.$ 
\end{minipage}
\hskip0.05\textwidth
\begin{minipage}[c]{.37\textwidth}
\begin{m-eqn}{
\scalebox{0.85}{
\begin{tikzpicture}[baseline=(baseline point)]
\coordinate (baseline point) at (0, 1.33);
\draw[thick,->] (-0.25,0) -- (4,0); \node at (3.85,0.25) {$\xx$};
\draw[thick,->] (0,-.05) -- (0,3.25); \node at (0.25,3) {$\yy$};
\node at (0,-0.25) {$(0,0)$};
\draw[densely dashdotdotted,thick] (3,0) -- (3,3);
\node at (1.25,1.5) {\scalebox{1.33}{$\yy = \ym$}};
\draw[thick,domain=0:3] plot (\x,{1/9*\x*\x*\x});
\draw[thick,densely dashdotdotted,domain=2.05:2.9] plot(\x,{3*\x-6.1});
\node at (2,-0.25) {$\PP$};\node at (3,-0.25) {$1$};\node at (3.55,3) {$(1,1)$};
\fill[pattern=dots] (3,0) -- (3,3) -- (2.05,0) -- cycle;
\end{tikzpicture}
}}\label{eq:B}
\captionof{figure}{\textscale{.9}{Geom.\,def.\,of \bl.}}\label{fig:tang}
\end{m-eqn}
\end{minipage}
\end{tabular}

Table~\ref{tab:bdry} lists values of $\yy_\exc$, for various $\nn, \kp=\nn^\aaa$ and $\xx$. (We use logarithmic scale for $\kp$.) The values $\aaa=3$ and $7$ were chosen to be about the same as the boundary layer exponent and twice as much. 

\renewcommand{\arraystretch}{1.15}{\textscale{.825}{
\begin{longtable}[H]
{|c|l|c|c|r|}
\cline{1-5}
$\nn$ 
&$\begin{array}{l}\;\aaa\,;\,\pp\\ (\kp\!=\!\nn^\aaa)\end{array}$  
&
{\begin{tabular}{c}$1-\xx=10^{-0.75{\pp}};$\\ $\yy_\exc(\xx)<\dots$\end{tabular}}
& 
{\begin{tabular}{c}$1-\PP=10^{-\pp};$\\ $\yy_\exc(\PP); \Delta(\yy)$\end{tabular}} 
&
{\begin{tabular}{c}$1-\xx=10^{-2\pp};$\\ $1-\yy_\exc(\xx);\;\Delta(\yy)$\end{tabular}}
\\ \hhline{=====}
\multirow{2}{*}{$\kern-1.5ex\begin{array}{c}4\\ \aph_\bl\!=\!3.4\end{array}\kern-1.5ex$}  
                & $3;\;1.67$ 
				& $5.5(-2);\;0.38$ 
				&$0.021;\;0.57;\;5(-2)$
 				&$4.4(-4);\;1.7(-2);\;3(-3)$
\\ \cline{2-5} 
				& $7;\;4.01$  
				& $9.7(-4);\;0.16$
				& $9.6(-5);\;0.54;\;1(-3)$  
				& $9.2(-9);\;9.5(-5);\;3.8(-7)$ 
\\ \hhline{=====} 
\multirow{2}{*}{$\kern-1.5ex\begin{array}{c}10\\ \aph_\bl\!=\!2.4\end{array}\kern-1.5ex$} 
                & $3;\;3.09$ 
				& $4.8(-3);\;0.60$ 
				& $8.1(-4);\;0.81;\;0.1$  
				&$6.5(-7);\;2.8(-4);\;5.3(-4)$
\\ \cline{2-5} 
				& $7$\;;\;$6.68$  
				& $9.6(-6);\;0.32$ 
				& $2.0(-7);\;0.70;\;2(-2)$  
				& $4.2(-14);\;1.8(-7);\;2.5(-8)$ 
\\ \hhline{=====} 
\multirow{2}{*}{$\kern-1.5ex\begin{array}{c}23\\ \aph_\bl\!=\!2.1\end{array}\kern-1.5ex$}
					& $3;\;5.49$ 
					& $7.6(-5);\;0.89$
					& $3.2(-6);\;0.99;\;0.2$  
					& $1.0(-11);\;3.6(-8);\;3.2(-6)$ 
\\ \cline{2-5}
					& $7;\;9.65$ 
					& $5.7(-8);\;0.56$
					& $2.1(-10);\;0.89;\;0.1$  
					& $4.8(-20);\;4.7(-11);\;1.7(-10)$
\\ \hhline{=====}
\multicolumn{1}{c}{Key:}& \multicolumn{4}{l}{for $(\nn,\aaa)=(4,7),$ $\;\kp=4^{7}=16384$, $\pp=4.01, \PP=1-10^{-4.01}\approx0.9999$} 
\\ \multicolumn{1}{l}{}& \multicolumn{4}{l}{At $1-\xx=10^{-2\pp}=9.2(-9)$, error is $\Delta(\yy)=\yip(\xx)-\yim(\xx)=3.8\cdot 10^{-7}$;} 
\\ \multicolumn{1}{l}{}& \multicolumn{4}{l}{$\yy_\exc(\xx)\approx1-9.5(-5)$ with $6$ certain decimals.}
\\ 
\caption{using envelopes for computations at the boundary layer}\label{tab:bdry}
\end{longtable}
}}\renewcommand{\arraystretch}{1}

The function $\yy_\exc$ drops steeply on the interval $[1-10^{-0.75\pp},1]$. The main information about the function is contained within the tiny interval $[\PP,1]$. 
Concerning the precision of the estimate $\yy_\exc(\xx)\in[\ym(\xx),\yp(\xx)]$, note that it is good for $\xx\ges 1-10^{-2\pp}$, but becomes loose for smaller $\xx$. Consequently, we divide $[0,1]$ in two pieces: for the computationally difficult boundary layer interval $[1-10^{-0.5\pp},1]$, we construct an explicit approximating formula for $\yy_\exc$; on $[0,1-10^{-0.5\pp}]$, $\yy_\exc$ has a gentle growth, so we can compute its values numerically (see~\S\ref{sssct:approx}). The envelopes are good at small $\xx$ only when $\nn$ is small compared to $\aaa$.  

\nit\begin{minipage}{0.97\textwidth}
\begin{minipage}[c]{.6\textwidth}
The maximal error $\yip-y_\exc$ over $[0,1]$ is $\yip(0)$. The smaller its value, the better the estimate on $[0,1]$. 

E.g. for $(\nn,\aaa)=(4,7)$, $\yy_\exc(0.95)\approx 0.01$, and $0.95$ is away from the boundary-interval $[0.9999,1]$.
\end{minipage}
\quad 
\begin{minipage}[c]{.37\textwidth}
\textscale{.825}{\renewcommand{\arraystretch}{1.15}
\begin{longtable}{|c|c|c|c|}
\cline{1-4} 
\diagbox[width=6ex, height=5ex]{$\aaa$}{$\nn$} & 4 & 7 & 10 
\\ \cline{1-4}
{7}  & 1(-3) & & 
\\ \cline{1-4}
{10} & 1(-4) & 1(-3) &  
\\ \cline{1-4}
{15} & 9(-7) & 5(-5) & 4(-4) 
\\ \cline{1-4}
{20} &9(-9) & 2(-6) & 3(-5) 
\\ \cline{1-4} 
\caption{error $\yip(0)$, small $\nn/\aaa$}\label{tab:y+}\label{tab:yOK}
\end{longtable} 
\renewcommand{\arraystretch}{1}
}
\end{minipage}
\end{minipage}

%%%%%%%%%%
%%%%%%%%%%

\subsubsection{Estimating \textit{C}\textsubscript{\rm ex}}

Its knowledge is essential, allows transforming the boundary- into an initial-value problem; for the latter there are robust algorithms. 
Let $\CC_-$ be as in Prop.~\ref{prop:estim-C}(iv) and $\ww_+$ be the solution of~\eqref{eq:h1}, with $\CC=\CC_-, \ww_+(1)=1$. Since $\CC_-\les\CC_\exc$, we have $\ww_+\ges \yy_\exc$; numerically, we use RK-method, we find that $\ww_+(0)>0$, in agreement with~\eqref{eq:yC}. 

\begin{m-corollary}\label{cor:C+}
The parameter $\CC_\exc$ in~\eqref{eq:h1} belongs to $[\CC_-,\Cp]$, so one has 
\\[.5ex]\centerline{
$\disp\CC_\exc=
\frac{1}{\nn+1}{\Bigl[ 1+\frac{\kp(\nn-1)}{\sqrt{2(\nn+1)}} \Bigr]}^{-\frac{2(\nn+1)}{\nn-1}+\gma_\exc},\; 0<\gma_\exc<\gma_+\!:=\!\frac{\ln((\nn+1)\Cp)}{\ln\Bigl( 1+\frac{\kp(\nn-1)}{\sqrt{2(\nn+1)}} \Bigr)}+\frac{2(\nn+1)}{\nn-1}.$
}
\end{m-corollary}
\nit For $(\nn,\aaa\!=\!\log_\nn(\kp))$, we compute with four decimals $\gma\!>\!0$ as follows. 

\nit\begin{tabular}{ll}
\begin{minipage}{.97\textwidth} 
\it\textscale{0.9}{The solution $\ww_\gma$ of   
$\;\left\{\begin{array}{l}
\disp\ww'=\sqrt{2}\kp\sqrt{\frac{\ww^{\nn+1}}{\nn+1}+\CC},\,\ww(1)\!=\!1,
\\ 
\disp\CC=\CC_\gma\!=\!\frac{1}{n+1}{\Bigl[ 1+\frac{\kp(\nn-1)}{\sqrt{2(\nn+1)}} \Bigr]}^{-\frac{2(\nn+1)}{\nn-1}+\gma}\!,
\end{array}\right.$ 
is so that: \\[.5ex] 
(i) $\ww_\gma(0)>0$;\; (ii) $\ww_{\gma'}(0)<0$ for $\gma'=\gma+10^{-4}$ (last decimal of $\gma$ increased by one, $\CC_\gma\mt\CC_{\gma'}$). 
\\(Note: for $\gma_{\blt\exc}$ with five/six decimals, one adds $10^{-5}/10^{-6}$.)
}
\end{minipage}
\end{tabular}\\[.5ex]
\nit In Table\ref{tab:ggw0}, $\gma$ with this property is denoted $\gma_{\blt\exc}$. Thus $\gma_\exc$ satisfies \\[.5ex]\centerline{$\gma_{\blt\exc}<\gma_\exc<\gma_{\blt\exc}+10^{-4}\;(\text{resp.}\,+10^{-5}/10^{-6})$.} 
\renewcommand{\arraystretch}{1.15}{\textscale{.825}{
\begin{longtable}[H]{c l|c| |l|c| |l|c| |l|c| |l|c|}
\cline{1-11}
  \multicolumn{1}{|c||}{\multirow{2}{*}{\diagbox[width=6ex, height=5ex]{$\aaa$}{$\nn$}}} &
  \multicolumn{2}{c||}{4} &
  \multicolumn{2}{c||}{7} &
  \multicolumn{2}{c||}{23} &
  \multicolumn{2}{c||}{73} &
  \multicolumn{2}{c|}{200} 
\\ \cline{2-11}
	\multicolumn{1}{|c||}{} 
& {$\kern0ex\aph_\bl\!=\!3.4\kern0ex$} & {$\ww(0)$} 
& {$\kern0ex\aph_\bl\!=\!2.5\kern0ex$} & {$\ww(0)$} 
& {$\kern0ex\aph_\bl\!=\!2.1\kern0ex$} & {$\ww(0)$} 
& {$\kern0ex\aph_\bl\!=\!2.3\kern0ex$} & {$\ww(0)$} 
& {$\kern0ex\aph_\bl\!=\!3.2\kern0ex$} & {$\ww(0)$} 
\\ \hhline{===========}
	\multicolumn{1}{|c||}{\multirow{2}{*}{0.5}} 
	& {2.6963} & {1(-6)} & {2.2301} & {1(-5)} 
	& {1.9303} & {2(-5)} & {1.8789} & {2(-4)} 
	& {1.87987} & {2(-5)} 
\\ 
	\multicolumn{1}{|c||}{} 
	& {3.1441} & {1.24} & {2.3898} & {1.08} 
	& {1.9573} & {0.97} & {1.8849} & {0.94} 
	& {1.88160} & {0.91} 
\\ \hhline{===========}
	\multicolumn{1}{|c||}{\multirow{2}{*}{1}} 
	& {1.7982} & {2(-6)} & {1.4520} & {3(-5)} 
	& {1.2476} & {2(-4)} & {1.2199} & {4(-5)} 
	& {1.22557} & {2(-5)} 
\\ 
	\multicolumn{1}{|c||}{} 
	& {2.0154} & {1.20} & {1.5248} & {0.99} 
	& {1.2599} & {0.88} & {1.2226} & {0.84} 
	& {1.22637} & {0.84} 
\\ \hhline{===========}
	\multicolumn{1}{|c||}{\multirow{2}{*}{3}}  
	& {0.6820} & {9(-6)} & {0.5676} & {5(-6)} 
	& {0.5084} & {3(-4)} & {0.5055} & {3(-4)} 
	& {0.51164} & {6(-5)} 
\\ 
	\multicolumn{1}{|c||}{} 
	& {0.7410} & {0.98} & {0.5919} & {0.99} 
	& {0.5133} & {1.06} & {0.5066} & {1.03} 
	& {0.51198} & {1.08} 
\\ \hhline{===========}
	\multicolumn{1}{|c||}{\multirow{2}{*}{7}} 
	& {0.2913} & {4(-7)} & {0.2527} & {2(-6)} 
	& {0.2324} & {1(-4)} & {0.2328} & {5(-5)} 
	& {0.23631} & {8(-5)} 
\\ 
	\multicolumn{1}{|c||}{} 
	& {0.3164} & {0.97} & {0.2635} & {1.03} 
	& {0.2346} & {1.11} & {0.2333} & {1.09} 
	& {0.23647} & {1.18} 
\\ \hhline{===========}
	\multicolumn{1}{|c||}{\multirow{2}{*}{15}} 
	& {0.1355} & {1(-9)} & {0.1197} & {1(-7)} 
	& {0.1114} & {2(-5)} & {0.11198} & {8(-6)} 
	& {0.113817} & {2(-5)} 
\\ 
	\multicolumn{1}{|c||}{} 
	& {0.1472} & {0.97} & {0.1249} & {1.06} 
	& {0.1125} & {1.19} & {0.11223} & {1.17} 
	& {0.113892} & {1.19} 
\\ \hhline{===========}
	\multicolumn{1}{|c||}{\multirow{2}{*}{20}} 
	& {0.1016} & {2(-12)} & {0.0901} & {3(-9)} 
	& {0.0840} & {1(-5)} & {0.08455} & {2(-5)} 
	& {0.085966} & {2(-5)} 
\\ 
	\multicolumn{1}{|c||}{} 
	& {0.1103} & {0.96} & {0.0940} & {1.06} 
	& {0.0848} & {1.15} & {0.08474} & {1.19} 
	& {0.086023} & {1.19} 
\\ \hhline{===========} 
\multicolumn{1}{r}{Key:} & \multicolumn{1}{|c|}{$\gma_{\blt\exc}\,$} & \multicolumn{3}{l|}{$\ww_{\gma_{\blt\exc}}(0)<\dots$} & \multicolumn{6}{r}{}
\\ \multicolumn{1}{r}{} &\multicolumn{1}{|c|}{$\gma_+$} & \multicolumn{3}{l|}{$(\gma_+-\gma_{\blt\exc})\cdot\aaa\cdot\nn\cdot\ln(\nn)$} & \multicolumn{6}{r}{}
\\ \cline{2-5} 
\caption{values of $\gma_{\blt\exc}, \gma_+$, $\ww_{\gma_{\blt\exc}}(0)=\mxerr$ on $[0,1]$.} 
\label{tab:ggw0}
\end{longtable}
}}\renewcommand{\arraystretch}{1}

\begin{m-remark} 
\nit(i) The graphs for the {\bvp}~\eqref{eq:h2} in Fig.~\ref{fig:ortshape} were plotted using the RK-method for the \ivp~\eqref{eq:h1}, with $\gma_{\blt\exc}$ in the table. They approach the origin, as required by $\yy(0)=0$. 

\nit(ii) The data shows that $\gma_+-\gma_{\blt\exc}$ decreases with $\nn,\kp$, thus $\yip$ approaches $y_\exc$.  
The visual explanation is that the graphs of $\yip, y_\exc$ become $\ort$-shaped, so the secant line defining $\Cp$ approaches the tangent to $\yy_\exc$ at the origin (see the geometric definition in Fig.~\ref{fig:Cp}). 

\nit(iii) If the precision given by the four/five-digit values of $\gma_{\blt\exc}$ is not satisfactory, one iterates shooting until the desired precision. Example: for $(\nn,\aaa)=(200,20)$, the value $0.085966$ yields $2.10^{-5}$-precision. Iteration gives $\gma_\exc\approx0.0859667086$ (very close to $\gma_+$) with $2.10^{-9}$-precision.
\end{m-remark}

\subsection{An approximating function}\label{ssct:appfct}

An explicit (accurate) approximation for the exact solution of~\eqref{eq:h2} is important for several reasons. First, it gives information about $\yy_\exc$, one can address further issues which may require explicit formulae. Second, it is essential from the perspective of the computational power and time. 

\subsubsection{An empirical estimate}\label{sssct:empir} 

Since $\gma_\exc$ approaches $\gma_+$ for large $\nn,\kp$, we estimated their difference (bottom-right entry of each cell in Table~\ref{tab:ggw0}). Surprisingly, these chaotic numbers obey a simple \emph{empirical} inequality: 
\\[.5ex]\centerline{
$0.5\les(\gma_+-\gma_{\blt\exc})\cdot\aaa\cdot\nn\ln(\nn)\les1.5
\;\;\Rightarrow\;\;\gma_+-\frac{2}{\aaa\nn\ln(\nn)}\les\gma_{\blt\exc}\les\gma_+.$ 
}
\begin{m-definition}\label{def:Cest}
Let $\gma_\gss:=\gma_+-\frac{1}{\aaa\nn\ln(\nn)}=\gma_+-\frac{1}{\nn\cdot\ln(\kp)}$ (average of bounds). 
We denote $\ww_{\gss}$  the numerical function solving~\eqref{eq:h1} with 
\\[.5ex]\centerline{
$\disp\CC=\CC_{\gss}
=\frac{1}{n+1}{\Bigl[ 1+\frac{\kp(\nn-1)}{\sqrt{2(\nn+1)}} \Bigr]}^{-\frac{2(\nn+1)}{\nn-1}+\gma_\gss}
=\Cp\cdot{\Bigl[ 1+\frac{\kp(\nn-1)}{\sqrt{2(\nn+1)}} \Bigr]}^{-\frac{1}{\nn\cdot\ln(\kp)}}.$
}\\[.5ex] 
As in Prop.~\ref{prop:estim-C}(iii), let $\disp\ell:=\frac{\yip(0)}{\yip(0)-\ww_{\gss}(0)}$, and consider the function: 
\begin{m-eqn}{
\YY_\gss:=(1-\ell)\yip+\ell \ww_{\gss}\;\; \text{(thus $\YY_\gss(0)=0,\,\YY_\gss(1)=1$).}
}\end{m-eqn} 
In most \bl-cases ($\aaa\ges3$) we have $\gma_\exc\les\gma_\gss$, thus  $\yy_\exc\ges\ww_{\gss}$;  Proposition~\ref{prop:estim-C} implies $\yy_\exc\ges\YY_\gss$. 
\end{m-definition}

\nit Below we computed $\YY_\gss(0.5)$ and $\yy_\exc(0.5)-\YY_\gss(0.5)$; the greatest error is attained about the middle of $[0,1]$, since $\YY_\gss, \yy_\exc$ coincide at $0, 1$. 
\renewcommand{\arraystretch}{1.15}{\textscale{.825}{
\begin{longtable}[H]{c l|c| |l|c| |l|c| |l|c| |l|c|}
\cline{1-11}
	\multicolumn{1}{|c||}{\diagbox[width=6ex, height=5ex]{$\aaa$}{$\nn$}} 
	&\multicolumn{2}{c||}{4} & \multicolumn{2}{c||}{7} 
	& \multicolumn{2}{c||}{23} 
	& \multicolumn{2}{c||}{73} 
	& \multicolumn{2}{c|}{200} 
\\ \hhline{===========}
	\multicolumn{1}{|c||}{\multirow{2}{*}{1}} 
	& {1.8350} & {0.38} & {1.4514} & {0.397} 
	& {1.2446} & {0.444} & {1.2194} & {0.474} 
	& {1.2254} & {0.488} 
\\ 
	\multicolumn{1}{|c||}{}
	& {1.7982} & {8(-4)} & {1.4520} & {2(-5)} 
	& {1.2476} & {8(-5)} & {1.2199} & {2(-5)} 
	& {1.2255} & {2(-4)} 
\\ \hhline{===========}
	\multicolumn{1}{|c||}{\multirow{2}{*}{3}} 
	& {0.6809} & {0.087} & {0.5674} & {0.122} 
	& {0.5086} & {0.255} & {0.5056} & {0.375} 
	& {0.51166} & {0.439} 
\\ 
	\multicolumn{1}{|c||}{}
	& {0.6820} & {7(-5)} & {0.5676} & {5(-6)} 
	& {0.5084} & {2(-4)} & {0.5055} & {2(-4)} 
	& {0.51164} & {4(-5)} 
\\ \hhline{===========}
	\multicolumn{1}{|c||}{\multirow{2}{*}{7}} 
	& {0.2906} & {2.23(-3)} & {0.2530} & {9.15(-3)} 
	& {0.2326} & {0.081} & {0.2328} & {0.232} 
	& {0.23633} & {0.355} 
\\ 
	\multicolumn{1}{|c||}{}
	& {0.2913} & {2(-7)} & {0.2527} & {3(-6)} 
	& {0.2324} & {7(-5)} & {0.2328} & {4(-5)} 
	& {0.23631} & {5(-5)} 
\\ \hhline{===========}
	\multicolumn{1}{|c||}{\multirow{2}{*}{15}} 
	& {0.1352} & {1.37(-6)} & {0.1200} & {5.10(-5)} 
	& {0.11158} & {8.34(-3)} & {0.11201} & {8.97(-2)} 
	& {0.113829} & {0.231} 
\\ 
	\multicolumn{1}{|c||}{}
	& {0.1355} & {4(-11)} & {0.1197} & {5(-8)} 
	& {0.11145} & {4(-6)} & {0.11198} & {2(-5)} 
	& {0.113817} & {2(-5)} 
\\ \hhline{===========}
	\multicolumn{1}{|c||}{\multirow{2}{*}{20}} 
	& {0.10135} & {1.35(-8)} & {0.09033} & {1.99(-6)} 
	& {0.08419} & {2.0(-3)} & {0.08458} & {4.94(-2)} 
	& {0.085975} & {0.177} 
\\ 
	\multicolumn{1}{|c||}{}
	& {0.10160} & {2(-12)} & {0.09014} & {8(-10)} 
	& {0.08409} & {6(-6)} & {0.08455} & {2(-5)} 
	& {0.085966} & {7(-5)} 
\\ \hhline{===========} 
\multicolumn{1}{r}{\textscale{1}{Key:}} & \multicolumn{1}{|c|}{\textscale{1}{$\gma_\gss$}} 
& \multicolumn{3}{l|}{\textscale{1}{$\YY_\gss(0.5)$}} &\multicolumn{6}{r}{} 
\\ 
\multicolumn{1}{r}{} & \multicolumn{1}{|c|}{\textscale{1}{$\gma_{\blt\exc}$}} 
& \multicolumn{3}{l|}{\textscale{1}{error=$|\yy_\exc(0.5)-\YY_{\gss}(0.5)|$}} &\multicolumn{6}{r}{} 
\\ \cline{2-5}
\caption{values of $\gma_\gss, \gma_{\blt\exc}$, $\YY_\gss(0.5)$, and the computational error.}\label{tab:gss}
\end{longtable}
}}\renewcommand{\arraystretch}{1}

\nit Most computations run within seconds. The data shows that the difference $\YY_\gss(0.5)-\yy_\exc(0.5)$ is (very) small, especially in the boundary layer situations, so $\yy_\exc$ is well-approximated by $\YY_\gss$. Overall, replacing $\yy_\exc$ with $\YY_\gss$, transforms the computationally challenging {\bvp}~\eqref{eq:h2} into a routine {\ivp}.

\medskip\begin{minipage}{.975\textwidth}\centering
\begin{minipage}[c]{.4\textwidth}
\centering\textscale{.8}{ 
{\includegraphics[height=0.62\textwidth]{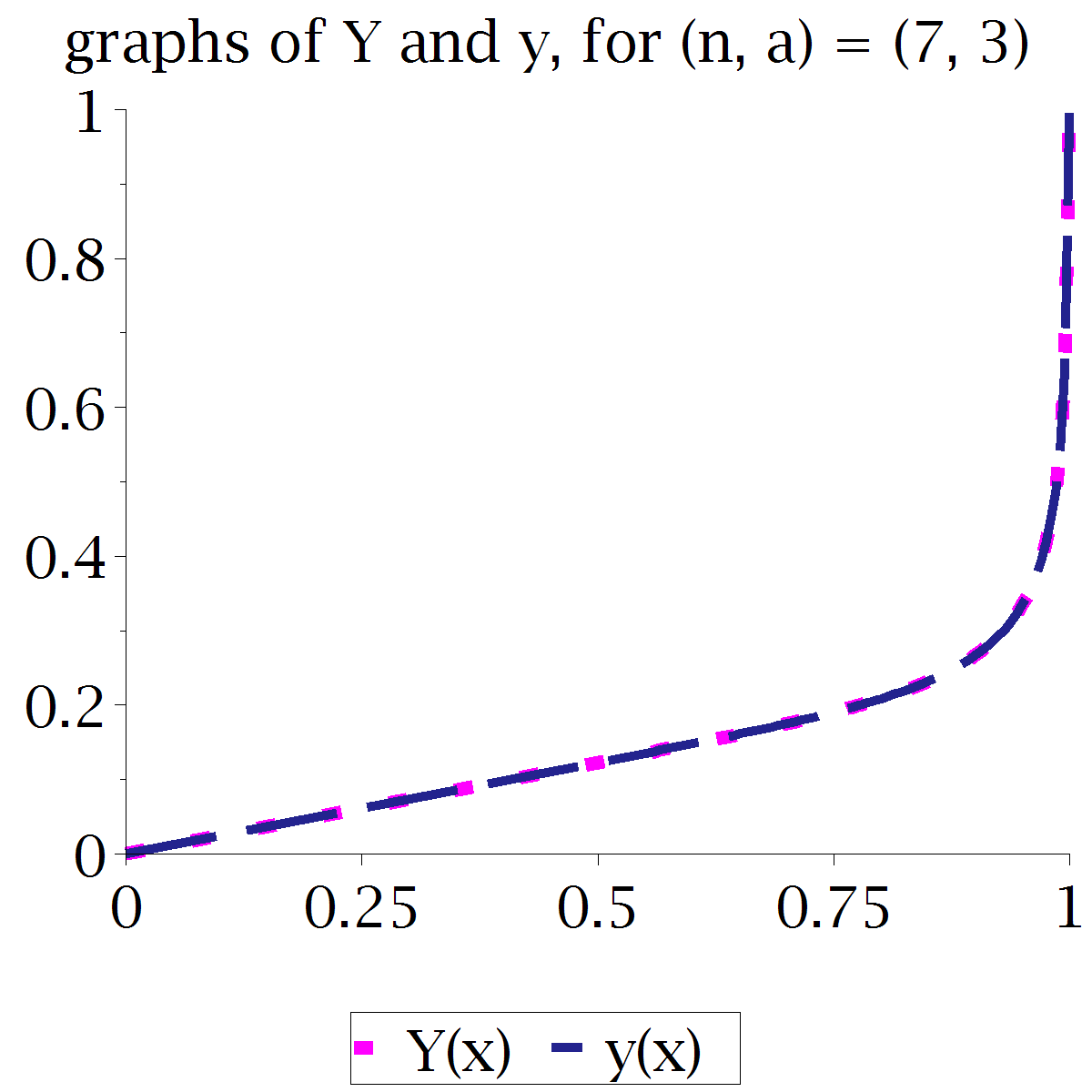}}
}
\end{minipage}
\begin{minipage}[c]{.4\textwidth}
\centering\textscale{.8}{ 
{\includegraphics[height=0.62\textwidth]{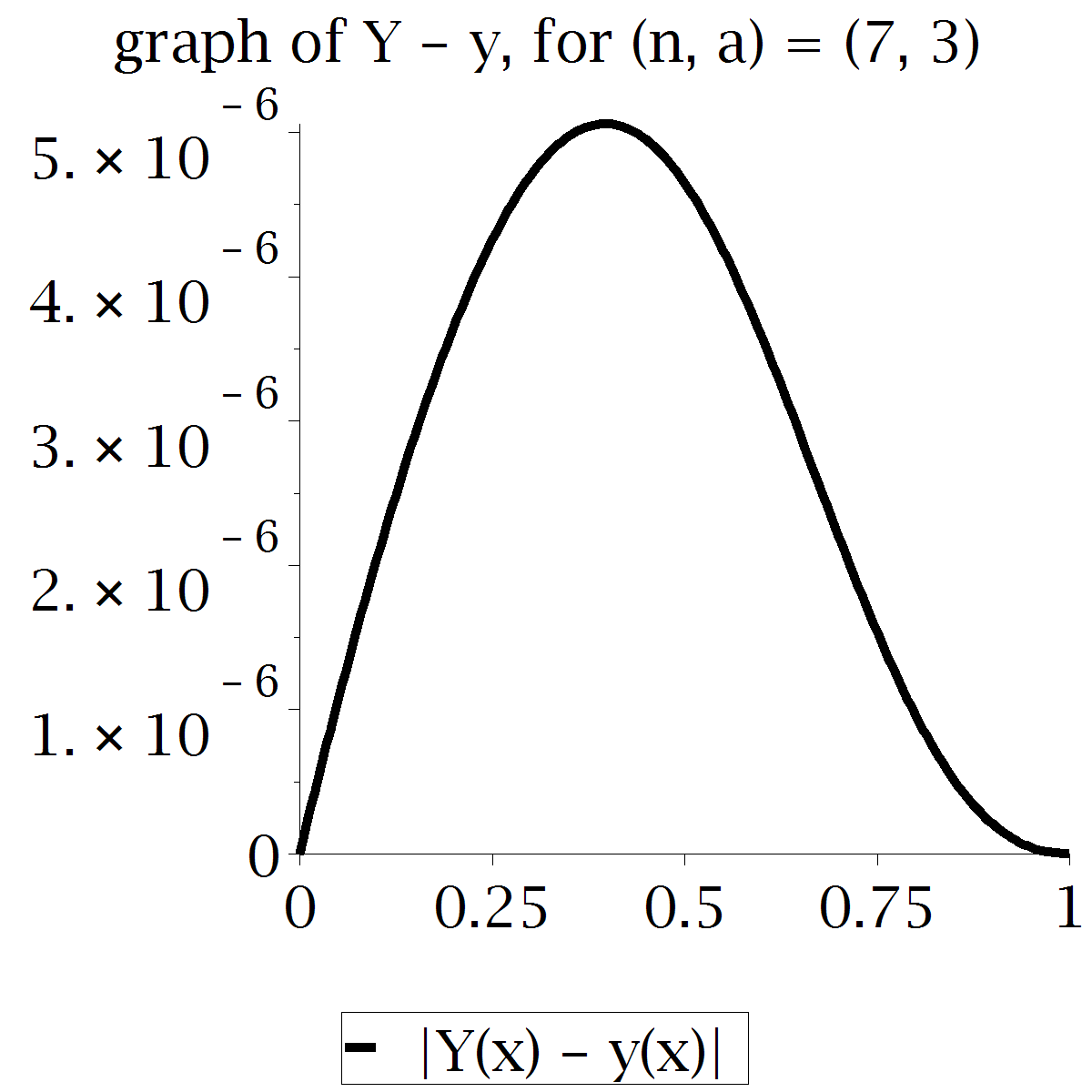}}
}
\end{minipage}
\captionof{figure}{graphs of $y_\exc$ and $\YY_\gss$, and $\YY_\gss-y_\exc$}\label{fig:Yy}
\end{minipage} 

\nit (For plotting the difference, we used high precision estimate for $\gma_{\blt\exc}$.)

%%%%%%%%%%
%%%%%%%%%%

\subsubsection{Defining the approximating function}\label{sssct:approx}

The explicit approximate value of $\CC_\exc$ allows defining a function $\ZZ$ which accurately approximates $\yy_\exc$ on an interval $[\aph_1,1]$ containing the boundary. Thus, to estimate $\yy_\exc(x)$ for some $x\in[0,1]$, we proceed as follows: 

{\rm(i)} for $\xx\in[\aph_1,1]$, simply evaluate $\ZZ(\xx)$;\;\; 
{\rm(ii)}  for $\xx\in[0,\aph_1]$, consider the {\ivp} 
\\[.75ex]\centerline{
\null\hfill
$\disp\,\disp\zz'=\frac{\kp\sqrt{2}}{\sqrt{\nn+1}}\cdot(\zz+c)^{\frac{\nn+1}{2}}, \;\zz(\aph_1)=\ZZ(\aph_1),$ 
\quad
\text{with}$\disp\;c=[1+(\nn+1)\CC_{\gss}]^{\frac{1}{\nn+1}}-1.$
\hfill$(\star)$
}\\[.75ex] 
Its starting point is at $x=\aph_1$ (not $x=1$). If $\ZZ_{num}$ denotes the numerical solution, $|\ZZ_{num}(0)|$ is the error on $[0,\aph_1]$. This approach bypasses the computationally sensitive boundary layer. 

\nit$\bullet$ The approximating function $\ZZ$ solves the {\ivp} $\disp\,\zz'=\frac{\kp\sqrt{2}}{\sqrt{\nn+1}}\cdot(\zz+c)^{\frac{\nn+1}{2}},\;\zz(1)=1,$ so 
\begin{m-eqn}{
\disp\ZZ(x):={\biggl[ (1+c)^{-\frac{\nn-1}{2}}+\frac{\kp(\nn-1)}{\sqrt{2(\nn+1)}}(1-\xx) \biggr]}^{-\frac{2}{\nn-1}}{\kern-1ex} -c.
}\label{eq:Z}
\end{m-eqn}
The value of $c$ in~$(\star)$ is obtained by imposing $\ZZ'(1)={\kp\sqrt{2}}\sqrt{\frac{1}{\nn+1}+\CC_{\gss}}$. With this choice, the right-hand side is almost $\yy'_\exc(1)={\kp\sqrt{2}}\sqrt{\frac{1}{\nn+1}+\CC_{\gma_\exc}}$. 
(Note that $\ZZ$ satisfies $\zz''=\kp^2(\zz+c)^\nn$, which is very close to $\yy''=\kp^2\yy^\nn$. The initial values of $\ZZ, \yy_\exc$ are almost equal, too. As a consequence, $\ZZ(x)-\yy_\exc(x)$ will be small for $x$ close to $1$.) 

We define $\aph_1:=1-10^{-0.5\pp}$ as in~\eqref{eq:B}, using the lower envelope of the exact solution. The reason for using the factor $0.5$ is to take $\alpha_1$ far out of the boundary layer region. Therefore, solving the {\ivp} $(\star)$ on $[0,\aph_1]$ creates no difficulties. 
Overall, we determined an explicit approximating function for the exact solution of the heat radiation equation. 

\medskip\nit$\bullet$ By iterating the construction, one defines the function $\ZZ$ on a longer interval. 
We consider the {\ivp} $(\star)$, but this time on the interval $[0,\aph_1]$ and with the initial value $\ZZ(\aph_1)$ at $x=\aph_1$. 

Note that this requires rescaling the corresponding {\bvp} to bring it into the form~\eqref{eq:h2}. The new constant $\kp_2=\kp\aph_1\ZZ(\aph_1)^{(\nn-1)/2}$ determines the new estimated value for $\CC$ (Definition~\ref{def:Cest}), hence of $c$ (in~$(\star)$). The lower envelope $\ym$ in Proposition~\ref{prop:estim-C} changes accordingly, and yields the point $\aph_2=\PP$, defined as in~\eqref{eq:B} (using the exponent $-\pp$ this time, for improved the accuracy, since we are outside the boundary layer). 

Briefly, the expression of $\ZZ$ on $[\aph_2,\aph_1]$ is the same~\eqref{eq:Z} but with modified $\kp$ and $c$. 

\begin{minipage}{.97\textwidth}
\begin{minipage}{.5\textwidth}
The concatenation of these two (explicit) expressions is a continuous, piecewise defined function which approximates $\yy_\exc$ on $[\aph_2,1]$. 

This two-step approximation largely suffices to cover the challenging bl-interval.
\end{minipage}
\hskip0.05\textwidth
\begin{minipage}{.4\textwidth}
\includegraphics[height=0.75\textwidth]{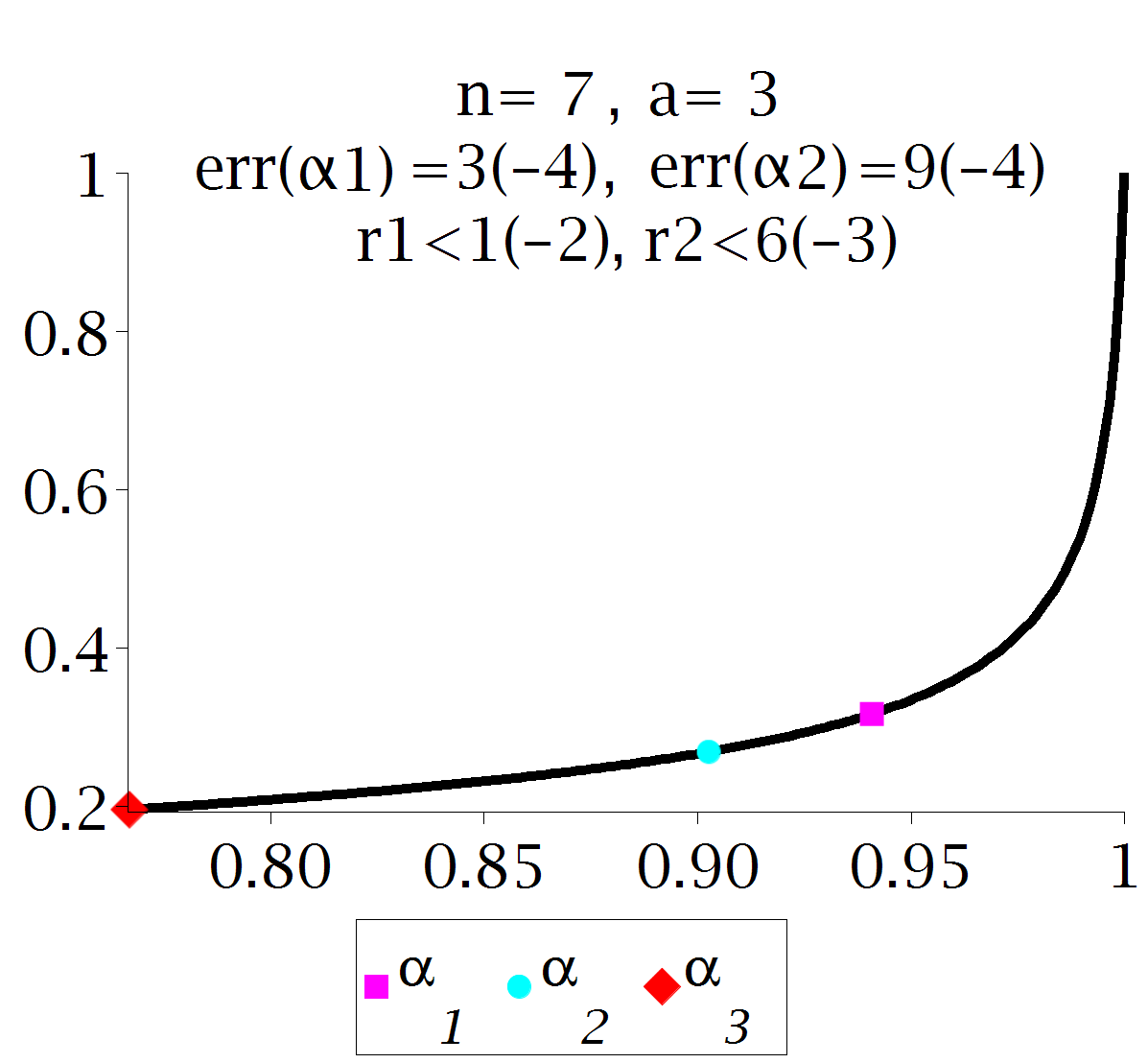}
\end{minipage}
\captionof{figure}{error $|\yy_\exc(\aph_i)-Z(\aph_i)|$ and ratio $r_i=\yy'_\exc(\aph_i)/\yy'_\exc(1)$, $i=1, 2$.}
\end{minipage}

\renewcommand{\arraystretch}{1.15}{\textscale{.825}{
\begin{longtable}{cr|rr}
\cline{1-4}
\multicolumn{1}{|c||}{\diagbox[width=6ex, height=5ex]{$\aaa$}{$\nn$}} 
&\multicolumn{1}{c|}{4}&\multicolumn{1}{c|}{23}&\multicolumn{1}{c|}{200}
\\ \hhline{====}%\cline{1-4}
\multicolumn{1}{|c||}{\multirow{2}{*}{3}}
&$0.218,\,6(-4)\;\&\;4(-4)$
&\multicolumn{1}{r|}{$0.679,\,3(-6)\;\&\;5(-4)$}
&\multicolumn{1}{r|}{$1-1.8\cdot10^{-11},\,2(-15)\;\&\;2(-4)$}
\\ \cline{2-4}
\multicolumn{1}{|c||}{}
&$0.149,\,2(-3)\;\&\;1(-3)$
&\multicolumn{1}{r|}{$0.677,\,3(-6)\;\&\;5(-4)$}
&\multicolumn{1}{r|}{$1-1.8\cdot10^{-11},\,2(-15)\;\&\;3(-4)$}
\\ \hhline{====}%\cline{1-4} 
\multicolumn{1}{|c||}{\multirow{2}{*}{7}}
&$0.034,\,2(-8)\;\&\;8(-6)$
&\multicolumn{1}{r|}{$0.336,\,5(-8)\;\&\;4(-4)$}
&\multicolumn{1}{r|}{$1-8.2\cdot10^{-5},\,2(-8)\;\&\;3(-4)$}
\\ \cline{2-4} 
\multicolumn{1}{|c||}{} 
&$0.002,\,4(-7)\;\&\;9(-6)$
&\multicolumn{1}{r|}{$0.332,\,7(-8)\;\&\;4(-4)$}
&\multicolumn{1}{r|}{$1-8.2\cdot10^{-5},\,2(-8)\;\&\;3(-4)$}
\\ \hhline{====}%\cline{1-4}
\multicolumn{1}{|c||}{\multirow{2}{*}{15}}
&$8.6\cdot10^{-4},\,3(-7)\;\&\;5(-8)$
&\multicolumn{1}{r|}{$0.101,\,4(-7)\;\&\;5(-5)$}
&\multicolumn{1}{r|}{$0.784,\,6(-8)\;\&\;3(-4)$}
\\ \cline{2-4} 
\multicolumn{1}{|c||}{} 
&$4.7\cdot10^{-4},\,2(-7)\;\&\;5(-8)$
&\multicolumn{1}{r|}{$0.090,\,2(-8)\;\&\;5(-5)$}
&\multicolumn{1}{r|}{$0.784,\,6(-8)\;\&\;3(-4)$}
\\\hhline{====}%\cline{1-4}
Key:&\multicolumn{2}{|c|}{$\ZZ(\aph_i),\,|\ZZ(\aph_i)-\yy_\exc(\aph_i)|\;\&\;|\ZZ_{num}(0)|$,\; for $i=1,2$}&
\\\cline{2-3}
\end{longtable}
\captionof{table}{values and error for the approximating function $\ZZ$ at $\aph_1$ (above) and $\aph_2$ (below).}\label{tab:ZZ}
}
We observe that, unless for $\nn\gg a$, the value of $\yy_\exc(\aph_2)$ is already low. For $\nn=4$ (case of the usual Stefan-Boltzmann equation), one reaches values less than $0.2$. 
For high values of $\nn$ (such as $200$), one may iterate the described process several times. 

%%%%%%%%%%
%%%%%%%%%%

\section{Conclusion}%\label{sct:concl}

We investigated the dependence on the parameters $\nn, \kp$ of the solution $\yy_\exc$ to the {\bvp} below, rooted in Stefan-Boltzmann law for radiative heat transfer: 
$$
\yy''=\kp^2\cdot \yy^\nn,\; \yy(0)=0,\; \yy(1)=1.
$$ 
To this end, we proceeded as follows.

\nit{\rm(i)}  
We set up a general analytical framework which allows constructing upper and lower envelopes to solutions of {\bvp}s. These are used to: 
\begin{itemize}
\item 
determine a priori upper/lower bounds for the initial derivative of the exact solution;
\item 
define the meaning of boundary layer for the exact solution;
\item 
decide what conditions lead to boundary layer.
\end{itemize}
\nit{\rm(ii)}  
We applied these techniques to the equation above. 
\begin{itemize}
\item 
Using explicit upper/lower envelopes, we defined and located the boundary-layer interval, and also computed $\yy_\exc$ at the boundary layer with high precision.
\item 
Determined a priori upper/lower bounds for $y_\exc'(0)$, in terms of $\kp,n$, which allow transforming the boundary- into an initial-value problem. The latter can be solved for a wide spectrum of parameters $\nn,\kp$ (on a regular laptop, with standard software). 
\item 
Numerical data shows that the (important, unknown) value of $\yy'_\exc(0)$ can be estimated by a simple analytical formula. This is used to construct an explicit function approximating $\yy_\exc$.
\end{itemize}

%\medskip\nit\textbf{Funding} This research did not receive funding.


\begin{thebibliography}{ooo}
\bibitem{amr} \bibauth{Ascher}{U.} \bibauth{Matheij}{R.} \bibauth{Russell}{R.} (1995) 
\bibbook{Numerical Solution of Boundary Value Problems for Ordinary Differential Equations}{SIAM, Classics in Appl. Math. 13}{Philadelphia}{2012} \url{https://doi.org/10.1137/1.9781611971231}

\bibitem{ch+pa} \bibauth{Choudhury}{S.} \bibauth{Paul}{R.} \bibtitl{A new approach to the generalization of Planck's law of black-body radiation} \bibjnyp{Ann. Phys.}{395}{2018}{317--325} \url{https://doi.org/10.1016/j.aop.2018.06.004}

\bibitem{gav} \bibauth{Gavalas}{G.} \bibbook{Nonlinear Differential Equations of Chemically Reacting Systems}{Springer Tracts Natural Phil. vol.17}{Springer-Verlag New York Inc.}{1968} \url{https://doi.org/10.1007/978-3-642-87643-1}

\bibitem{hal-aip} \bibauth{Halic}{M.} \bibtitl{An Approximate Solution to the One-Dimensional Heat Conduction Equation With Radiation} \bibjnyp{AIP Conf. Proc.}{3315}{2025}{400017} \url{https://doi.org/10.1063/5.0286522}

\bibitem{keg} \bibauth{Khosropour}{B.} \bibauth{Eghbali}{M.} \bibauth{Ghorbanali}{S.} \bibtitl{q-nonlinear Schrodinger and q-nonlinear Klein–Gordon equations in the frame work of GUP} \bibjnyp{Gen. Relativ. Gravity}{50}{2018}{25} \url{https://doi.org/10.1007/s10714-018-2348-6}

\bibitem{len+men} \bibauth{Lenzi}{E.} \bibauth{Mendes}{R.} \bibtitl{Blackbody radiation in nonextensive Tsallis statistics: Exact solution} \bibjnyp{Phys. Lett.~A}{250}{1998}{270--274} \url{https://doi.org/10.1016/S0375-9601(98)00819-6}

\bibitem{lvc} \bibauth{Lucchesi}{C.} \bibauth{Vaillon}{R.} \bibauth{Chapuis}{P.-O.} \bibtitl{Temperature dependence of near-field radiative heat transfer above room temperature} \bibjnyp{Materials Today Phys.}{21}{2021}{100562} \url{https://doi.org/10.1016/j.mtphys.2021.100562}

\bibitem{mppt} \bibauth{Martinez}{S.} \bibauth{Pennini}{F.} \bibauth{Plastino}{A.} \bibauth{Tessone}{C.} \bibtitl{q-Thermostatistics and the black-body radiation problem} \bibjnyp{Physica~A}{309}{2002}{85--105} \url{https://doi.org/10.1016/S0378-4371(02)00621-0}

\bibitem{me-ar} \bibauth{Mehta}{B.} \bibauth{Aris}{R.} \bibtitl{A Note on a Form of the Emden-Fowler Equation} \bibjnyp{J.~Math. Anal. Appl.}{36}{1971}{611-621} \url{https://doi.org/10.1016/0022-247X(71)90043-6}

\bibitem{mod} \bibauth{Mazumder}{S.} \bibauth{Modest}{M.} \bibbook{Radiative heat transfer, $4^{\rm th}\;{\rm ed.}$}{Acad. Press}{London}{2023} \url{https://doi.org/10.1016/C2018-0-03206-5}

\bibitem{nrt} \bibauth{Nobre}{F.} \bibauth{Rego-Monteiro}{M.} \bibauth{Tsallis}{C.} \bibtitl{Nonlinear Relativistic and Quantum Equations with a Common Type of Solution} \bibjnyp{Phys. Review Lett.}{106}{2011}{140601} \url{https://doi.org/10.1103/PhysRevLett.106.140601}

\bibitem{zzw} \bibauth{Zhu}{T.} \bibauth{Zhang}{Y.-M.} \bibauth{Wang}{J.-S.} \bibtitl{Super-Planckian radiative heat transfer between coplanar two-dimensional metals} \bibjnyp{Phys. Rev. B}{109}{2024}{245427} \url{https://doi.org/10.1103/PhysRevB.109.245427}
\end{thebibliography}
\end{document}